# On the Connection Between Irrationality Measures and Polynomial Continued Fractions

Nadav Ben David[1], Guy Nimri[1], Uri Mendlovic[2], Yahel Manor[1], Carlos De la Cruz Mengual[1], and Ido Kaminer[1]

[1]Technion - Israel Institute of Technology, Haifa 3200003, Israel
[2]Google Inc., Tel Aviv 6789141, Israel

### Abstract

Linear recursions with integer coefficients, such as the recursion that generates the Fibonacci sequence $F_n = F_{n-1} + F_{n-2}$, have been intensely studied over millennia and yet still hide interesting undiscovered mathematics. Such a recursion was used by Apéry in his proof of the irrationality of $\zeta(3)$, which was later named the Apéry constant. Apéry's proof used a specific linear recursion that contained integer polynomials (polynomially recursive) and formed a continued fraction; such formulas are called polynomial continued fractions (PCFs). Similar polynomial recursions can be used to prove the irrationality of other fundamental constants such as $\pi$ and $e$. More generally, the sequences generated by polynomial recursions form Diophantine approximations, which are ubiquitous in different areas of mathematics such as number theory and combinatorics. However, in general it is not known which polynomial recursions create useful Diophantine approximations and under what conditions they can be used to prove irrationality.

Here, we present general conclusions and conjectures about Diophantine approximations created from polynomial recursions. Specifically, we generalize Apéry's work from his particular choice of PCF to any general PCF, finding the conditions under which a PCF can be used to prove irrationality or to provide an efficient Diophantine approximation. To provide concrete examples, we apply our findings to PCFs found by the Ramanujan Machine algorithms to represent fundamental constants such as $\pi$, $e$, $\zeta(3)$, and the Catalan constant. For each such PCF, we demonstrate the extraction of its convergence rate and efficiency, as well as the bound it provides for the irrationality measure of the fundamental constant. We further propose new conjectures about Diophantine approximations based on PCFs.

Looking forward, our findings could motivate a search for a wider theory on sequences created by any linear recursions with integer coefficients. Such results can help the development of systematic algorithms for finding Diophantine approximations of fundamental constants. Consequently, our study may contribute to ongoing efforts to answer open questions such as the proof of the irrationality of the Catalan constant or of values of the Riemann zeta function (e.g., $\zeta(5)$).

# 1 Introduction

## 1.1 Apéry’s constant and his polynomial continued fraction

In his paper [1,2], R. Apéry ingeniously presented a specific linear recursion with integer polynomial coefficients which is used to prove the irrationality of $\zeta(3)$. This polynomial recursion generated two sequences $p_n, q_n$ (given different initial values) such that $p_n/\,q_n \underset{n\to\infty}{\longrightarrow} 6/\zeta(3)$, i.e., which constituted a Diophantine approximation of $\zeta(3)$. Apéry then showed that this specific sequence $p_n/\,q_n$ proved the irrationality of the number to which it converges. It is also demonstrated [2] that the linear recursion is equivalent to the following polynomial continued fraction (PCF):

$$\frac{6}{\zeta(3)} = 5 - \cfrac{1}{117 - \cfrac{64}{535 - \cfrac{729}{1463 \ldots - \cfrac{n^6}{34n^3 + 51n^2 + 27n + 5 \ldots}}}}.$$

Apéry’s paper inspired other researchers to apply related strategies to other problems in Diophantine approximations, to study irrationality measures of other constants, and to find applications in other fields [3-10].

Apéry’s result hints at a more general question: Which PCFs prove the irrationality of the number to which they converge? In other words: Which pairs of integer polynomials (such as $-n^6$ and $34n^3 + 51n^2 + 27n + 5$ in Apéry’s case) can be used to prove irrationality? This question is directly related to the intrinsic properties of PCFs, specifically, their rate of convergence and the properties of the Diophantine approximation sequences they create.

## 1.2 Polynomial continued fractions

In their most general form, PCFs are generalized continued fractions whose terms $a_n = a(n)$ and $b_n = b(n)$ are defined by two polynomials with integer coefficients, $a(n)$ and $b(n)$:

$$\mathrm{PCF}[a_n, b_n] = a_0 + \cfrac{b_1}{a_1 + \cfrac{b_2}{a_2 + \cfrac{b_3}{a_3 \ldots + \cfrac{b_n}{a_n \ldots}}}}$$

$$a_n = a(n), b_n = b(n) \in \mathbb{Z}[n].$$

Truncating a PCF at a finite step $n$ produces a rational number

$$\frac{p_n}{q_n} = a_0 + \cfrac{b_1}{a_1 + \cfrac{b_2}{a_2 + \cfrac{b_3}{a_3 \ldots + \cfrac{b_n}{a_n}}}},$$

called the $n$-th convergent of the PCF. The $n$-th convergent also results as the quotient of the numbers $p_n$ and $q_n$ defined independently as the solutions of the following linear recursion of depth 2

$$u_n = a_n u_{n-1} + b_n u_{n-2}, \qquad \text{(Eq. 1)}$$

with initial conditions,

$$p_{-1} = 1, \;\; p_0 = a_0$$

$$q_{-1} = 0, \;\; q_0 = 1.$$

We say that the value of the PCF is the limit of the sequence $p_n/q_n$, if it exists. There exist Mobius transformations with integer coefficients that transform between the limits of $p_n/q_n$ for different initial conditions – for any two pairs of rational, linearly-independent initial values [11].

PCFs appear in a wide range of fields of mathematics and are related to many special functions, including all trigonometric functions, exponentials, Bessel functions, generalized hypergeometric functions, the Riemann zeta function, and many other important functions such as erf and log [12-14]. Moreover, any infinite sum can be converted into a continued fraction, known as Euler's continued fraction. The space of PCFs also contains all linear recursions of depth 2 with rational polynomial coefficients (and some of their generalizations). In his study, Apéry developed a linear recursion of depth 2 with rational polynomials, which can be converted to a PCF, using the standard definition above.

### 1.3 The goals of this paper

Looking at the bigger picture, it is interesting to generalize Apéry's PCF. Consider an arbitrary linear integer recursion (of any order) used to create the numerators and denominators in a sequence of rational numbers. In other words, provided two sets of initial conditions, for the numerator and denominator, the linear recursion creates a Diophantine approximation sequence. Each such sequence may provide an efficient representation of the limit of the sequence. Intuitively, the efficiency is described by the rate (as a function of $n$) at which the sequence converges relatively to the sizes of the denominators. What can be said about the resulting sequence? What condition should the linear recursion fulfill for the generated sequence to prove that its limit is irrational? More generally, what bounds on irrationality measures does each linear recursion create?

In this paper, we describe the construction of a systematic method to find, for each PCF, the efficiency of its limit approximation, i.e., the lower bound it provides for the irrationality measure (we address a lower bound that simultaneously provides an upper bound [2]). We develop a criterion on the PCF for proving the irrationality of its limit. Specifically, **Theorem 2** states a formula for the irrationality bound for each PCF, yielding $\frac{ln|\alpha|-ln|B|+ln\,\lambda}{ln|\alpha|-ln\,\lambda}$, where $\alpha$ and $B$ can be calculated directly from $a_n$'s and $b_n$'s coefficients and $\lambda$ relates only to the growth rate of the PCF's greatest common divisor, denoted as

$$\mathrm{GCD}_n \stackrel{\mathrm{def}}{=} \mathrm{GCD}[p_n, q_n]$$

(specifically, $\ln(\lambda) = \limsup_{n\to\infty} \frac{1}{n} \ln\left(\frac{\mathrm{GCD}_n}{n!^{\deg a_n}}\right)$). A major advantage of this formula is that it does not require the PCF's limit as its input. Moreover, **Conjecture 1.1** states that, for the growth rate of $\mathrm{GCD}_n$ to be sufficient for an irrationality proof, the polynomial $b_n$ must have rational roots only.

PCFs that yield efficient Diophantine approximations are in general also better for computing more quickly the numbers to which they converge. Consequently, the results of our study could be used to develop faster means for high precision calculations of fundamental constants, such as attempts to compute more digits and study the normality of such constants [15-22].

Any mathematical expression that can be converted to PCFs, such as infinite sums used for the computation of fundamental constants [15-17,20,21], could be analyzed with the approach that we present in this paper. The conjectures that arise from our study hint at a general theory that goes beyond PCFs to any polynomial recursion, and maybe eventually beyond it to any linear recursion with rational coefficients.

Some of the conclusions of our study presented below go beyond PCFs and the motivation of irrationality proofs. In general, when given any linear recursive formula with integer coefficients, not necessarily one representing a PCF, it is interesting to study the greatest common divisor of two (or maybe more) sequences arising from the same recursion with different initial conditions. We find the solution for special cases of linear recursions, showing the rate of growth of $\mathrm{GCD}_n$. We hope that our study will contribute to efforts toward finding the general rules for greatest common devisors of arbitrary linear recursions.

### 1.4 Motivation and potential applications

Many of the PCF formulas that led us to the conjectures and proofs in this paper were originally found by the Ramanujan Machine project [10,30], which employs computer algorithms to find conjectured formulas for fundamental constants. Various algorithms are being developed as part of that project, and so far, they all focus on formulas in the form of PCFs. Since the algorithms check candidate formulas by their numerical match to target constants, the results are always in the form of conjectures rather than proven theorems. The first algorithms succeeded in finding conjectured PCF formulas for $\pi$, $e$, values of the Riemann zeta function $\zeta$, and the Catalan constant [10]. These latter formulas led to a new world record for the irrationality bound of the Catalan constant. The theorems and conjectures below can also help improve future algorithms that search for such conjectures.

We point out three interesting challenges that motivate this work, each having prospects in Diophantine approximations, as well as in experimental mathematics, i.e., computation-driven mathematical research (e.g., [10,23,24]):

(1) Given the polynomials $a_n, b_n$, determine whether $\mathrm{PCF}[a_n, b_n]$ provides a bound on the irrationality measure of its limit (if it converges), and if so, find the bound analytically from $a_n, b_n$.

(2) Estimate the efficiency of each PCF for computing fundamental constants to high precision.

(3) Develop faster algorithms to compute PCFs; more generally, compute any polynomial recursion more efficiently.

### 1.5 The measure of irrationality

The irrationality measure of a number $L$ is the largest $\delta$ for which there exists an infinite sequence of rational numbers $p_n/q_n \neq L$ such that

$$\left|L - \frac{p_n}{q_n}\right| < \frac{1}{|q_n|^{1+\delta}}. \qquad (\mathrm{Eq.}\,2)$$

This maximal $\delta$ is called the irrationality measure of the number $L$ [2,15], or the Liouville–Roth exponent. For irrational numbers, this maximum can be obtained by the *regular* continued fraction of $L$; however, its closed formula is often unknown (e.g., in the case of $\pi$). The Diophantine approximation is thought of as more *efficient* when $\delta$ is larger. Rational numbers have an irrationality measure 0, meaning that they cannot be approximated efficiently by other rational numbers. This property is part of the irrationality criterion: if there exists a sequence $p_n/q_n$ for which the inequality in $(\mathrm{Eq.}\,2)$ holds for some $\delta > 0$, then $L$ is irrational. Further, if there exists a sequence $p_n/q_n$ for which this inequality holds for some $\delta > 1$, then $L$ is transcendental by the Siegel–Roth theorem. Finally, if the inequality holds for arbitrarily large values of $\delta$, then $L$ is a Liouville number (infinite irrationality measure) [15]. Intuitively, for a sequence that proves irrationality, the growth of the denominator should be sufficiently slow in relation to the convergence rate of the PCF.

For each number $L$, there always exists a sequence that satisfy $(\text{Eq.}\,2)$ with $\delta = 0$ (for a rational $L$) or $\delta \geq 1$ (for an irrational $L$). However, the largest *known* $\delta$ can be smaller or larger than $0$. To find even one explicit sequence that reaches the maximal value is challenging. This challenge continues to motivate searches for new sequences $p_n/q_n$, from which one can extract larger lower bounds for the irrationality measures of constants. Each constant for which the rationality or irrationality is still unknown has all its known sequences with $\delta \leq 0$ (as in the case of the Catalan constant [19-22]). Then, finding one sequence for which $\delta > 0$ will directly prove irrationality. When $\delta$ is known to be positive, as in $\pi$, it is still interesting to find better PCFs with larger $\delta$ values, because it improves the bounds on the constant's irrationality measure (e.g., $\pi$'s upper bound [25] by Zeilberger and Zudilin). Therefore, it is of interest to find sequences for which $\delta$ is as large as possible, even when the value is negative.

In the rest of this paper, we use the symbol $\delta = \delta_{p_n,q_n}$ to denote the lower bounds on irrationality measures that we determine by analyzing sequences $p_n/q_n$ generated by PCFs. We identify a necessary and sufficient condition for the growth rate of $\text{GCD}_n$, under which $\text{PCF}[a_n, b_n]$ provides a nontrivial $\delta$ (larger than $-1$). Utilizing this criterion, we present an expression for $\delta$ and propose conjectures for its dependence on the choice of $a_n, b_n$.

# 2 Results

## 2.1 Summary of the main results

Unless stated otherwise, we focus on "balanced-degree" PCFs, defined as those for which $\frac{\deg b_n}{\deg a_n} = 2$. This PCF type is arguably the most common in the literature related to mathematical constants (see Appendix A, Ref. [10], and further references therein). We show that the growth rate (as a function of $n$) of $\mathrm{GCD}_n$ is key to the analysis of PCFs of this type. We find special interest in cases of $\mathrm{PCF}[a_n, b_n]$ for which $\mathrm{GCD}_n$ upholds $\limsup_{n\to\infty} \frac{1}{n} \ln\left(\frac{\mathrm{GCD}_n}{n!^{\deg a_n}}\right) > -\infty$. That is, while the denominator $q_n$ can be shown to grow super-exponentially as $n!^{\deg a_n}$, the reduced denominator $\frac{q_n}{\mathrm{GCD}_n}$ is of exponential order only. We call this phenomenon **factorial reduction (FR)**.

Below, we prove that for a PCF to provide a nontrivial $\delta$ value, FR is necessary (**Theorem 1**). We also derive formulas for this $\delta$ (**Theorem 2**), which could help provide irrationality proofs. The other results of our work are conjectures which attempt to provide a complete characterization of PCFs with nontrivial $\delta$s. All the conjectures are backed with extensive, computer-based, numerical tests and await formal proof.

Numerical tests show that FR is possible if and only if $b_n$ has rational roots only (**Conjecture 1.1**). We further conjecture that for every such $b_n$, the $\mathrm{PCF}[a_n, b_n]$ has FR for infinitely many choices of polynomials $a_n$ with *rational* coefficients (**Conjecture 1.2**). Other experiments revealed that, for the special case of $\deg b_n = 2$ and $\deg a_n = 1$, each $b_n$ with all-rational roots has exactly two infinite families of integer polynomials $a_n$ for which $\mathrm{PCF}[a_n, b_n]$ has FR (**Conjecture 1.4** presents the formula for the $a_n, b_n$ pairs). A necessary and sufficient condition on $a_n, b_n$ of arbitrary degrees for PCF$[a_n, b_n]$ to have FR still awaits discovery and a proof.

Example: Apéry's PCF

Let us explain how Apéry's PCF appears as a special case in our study. First, his PCF uses $b_n = -n^6$ which has only rational roots (a special case of **Conjecture 1.1**). Second, Apéry's PCF has FR and provides a nontrivial $\delta$ (a special case of **Theorem 1**), as he proved that a lower bound for the size of $\text{GCD}_n$ is $n!^3/e^{3n}$. Third, Apéry proved that $\delta = \frac{\ln \alpha - 3}{\ln \alpha + 3}$, for $\alpha = 17 + 12\sqrt{2}$, which exactly matches our general expression (see **Theorem 2**). Below, we generalize this process and conclusions to all PCFs, classify their different $\text{GCD}_n$s, and present a criterion for the PCF that allows us to prove irrationality.

## 2.2 Theorems about factorial reduction

We tested many PCFs for FR and identified a surprising phenomenon: despite the rarity of FR in an experimentally random PCF, we have so far found FR in *every* PCF that converges to a *fundamental constant* we tested (PCFs that converge to $\pi$, $e$, $\zeta(3)$, $\zeta(5)$, and the Catalan constant). Specifically, we tested *all* the PCFs found so far in the Ramanujan Machine project [10,30] and many other PCF formulas. This relation between FR and PCFs of fundamental constants is surprising because the algorithmic search in [10] did not favor PCFs that have FR. This intriguing fact hints at an underlying structure of PCFs that is required for formulas that converge to certain mathematical constants.

**Theorem 1** *(The necessity of FR). Consider a balanced-degree* $PCF[a_n, b_n]$ *satisfying* $A^2 + 4B > 0$. *For the* PCF *to provide a nontrivial* $\delta$*, it must have FR; i.e.,*

$$\limsup_{n\to\infty}\frac{1}{n}\ln\left(\frac{\mathrm{GCD}_n}{n!^{\deg a_n}}\right) > -\infty.$$

*In other words,* $\mathrm{GCD}_n$ *divided by* $n!^{\deg a_n}$ *is of exponential order, and so, the reduced denominators* $\frac{q_n}{\mathrm{GCD}_n}$ *are of exponential order rather than factorial.*

Satisfing $A^2 + 4B > 0$ ensures the PCF's convergence (see Appendices A and B).

*Outline of the proof.*

An outline of the proof is presented here. For the complete proof See Appendix B.

1. Bound the finite calculation error $\left|\frac{p_n}{q_n} - \lim_{n\to\infty}\frac{p_n}{q_n}\right|$ by $\left|\frac{\prod_{i=1}^{n} b_i}{q_{n+1}q_n}\right|$ (**Lemma 2**)
2. Extract an expression for $\delta$ from $\left|\frac{\prod_{i=1}^{n} b_i}{q_{n+1}q_n}\right| < \frac{1}{\left|\frac{q_n}{\mathrm{GCD}_n}\right|^{1+\delta}}$ (**Lemma 1**)
3. Estimate $\left|\prod_{i=1}^{n} b_i\right|$ (**Lemmas 3**) and $|q_{n+1}q_n|$ (**Lemmas 4**, using Poincaré's [31])
4. Conclude that if $\mathrm{GCD}_n$ does not grow super-exponentially as $n!^{\deg a_n}$, then $\delta = -1$

□

To size $\mathrm{GCD}_n$ for *any* PCF, we define $\lambda$ as the exponential order of $\mathrm{GCD}_n$ divided by the factorial part, i.e. (denoting $\ln 0 = -\infty$),

$$\ln(\lambda) \stackrel{\text{def}}{=} \limsup_{n\to\infty}\frac{1}{n}\ln\left(\frac{\mathrm{GCD}_n}{n!^{\deg a_n}}\right).$$

Hence, FR in fact means that $\lambda > 0$. So, assuming the above limit converges, we can write

$$\mathrm{GCD}_n \doteq \lambda^n \cdot n!^{\deg a_n}.$$

The notation $x_n \doteq y_n$ represents that for some $\kappa \in \{-1,1\}$

$$\lim_{n\to\infty}\frac{1}{n}\ln\kappa\frac{x_n}{y_n}=0. \qquad \text{(Eq. 3)}$$

This is an equivalence relation implying that $x_n$ and $y_n$ exhibit identical exponential growth rates as $n$ approaches infinity, although they may differ by sign, and factors that grow (or decay) more slowly than exponentials, such as polynomial factors ($n^r \doteq 1$ for any $r \in \mathbb{R}$). This relation is particularly useful for comparing sequences in our analysis, where the exponential terms are the dominant component and sub-exponential terms are negligible. Also, note that we use the fact that $x_n \doteq y_n$ and $x'_n \doteq y'_n$ implies $x_n x'_n \doteq y_n y'_n$.

**Theorem 2** *(A formula for $\delta$). Consider a balanced-degree* $\mathrm{PCF}[a_n, b_n]$ *with* $A^2+4B>0$ *converging to a number L* $\in \mathbb{R}$*. The irrationality measure of L is at least*

$$\delta=\frac{\ln|\alpha|-\ln|B|+\ln\lambda}{\ln|\alpha|-\ln\lambda},$$

*where A and B are the leading coefficients of* $a_n$ *and* $b_n$*, respectively, and* $\alpha$ *is the larger solution in absolute value of the equation*

$$x^2-Ax-B=0.$$

This formula establishes a link between the size of $\mathrm{GCD}_n$ and the value of $\delta$. As anticipated, the larger $\mathrm{GCD}_n$ is, the greater $\delta$ becomes, i.e., larger values of $\lambda$ imply larger values of $\delta$, which is advantageous for proving irrationality. The minimal possible $\lambda$ value is $\lambda=0$ which implies no FR and a trivial $\delta=-1$, while the maximal value is $\lambda=|\alpha|$ (as Appendix B shows $p_n \doteq q_n \doteq \alpha^n \cdot n!^{\deg a_n}$), which implies an infinite irrationality measure, i.e., a Liouville number. Any $\lambda$ value exhibiting $\lambda > \left|\frac{B}{\alpha}\right|$ implies $\delta>0$ and proves the irrationality of the PCF's limit.

*Outline of the proof.* Following **Theorem 1**, the formula for $\delta$ is derived from assuming $\mathrm{GCD}_n \doteq \lambda^n \cdot n!^{\deg a_n}$ and inserting it into the expression from **Lemma 1** (step 2 at the outline of **Theorem 1**'s proof). For the complete proof See Appendix B. □

Examples: Different values of $\lambda$ and $\delta$

1. Apéry's: $b_n = -n^6$ and $a_n = 34n^3 + 51n^2 + 27n + 5$, and therefore, $\alpha = \frac{34+\sqrt{34^2-4}}{2}$. Apéry showed [1] that $\left(\frac{n!}{\mathrm{LCM}[n]}\right)^3 | \mathrm{GCD}_n$ and thus $\mathrm{GCD}_n \geq \left(\frac{n!}{\mathrm{LCM}[n]}\right)^3$, where $\mathrm{LCM}[n]$ is the least common multiple (LCM) of $1,2 \dots n$, which satisfies $\mathrm{LCM}[n] \doteq e^n$ (as it is equal to the exponent of the second Chebyshev function [2]). Therefore, $\lambda \geq \frac{1}{e^3}$ and our formula provides $\delta = \frac{\ln|\alpha| - \ln|B| + \ln\lambda}{\ln|\alpha| - \ln\lambda} \geq 0.080529 \dots$, which exactly matches the Apéry's result.
2. Other irrational limits: For any integer $k \geq 3$, take $b_n = -n^2$ and $a_n = k(2n+1)$; therefore, $\alpha = k + \sqrt{k^2-1}$. The proof in Appendix D shows that $\frac{n!}{\mathrm{LCM}[n]} | (\mathrm{GCD}_n \cdot 2^n)$ and thus $\mathrm{GCD}_n \geq \frac{n!}{2^n \mathrm{LCM}[n]}$. Therefore, $\lambda \geq \frac{1}{2e}$, providing $\delta \geq \frac{\ln\left(k+\sqrt{k^2-1}\right) - \ln 2 - 1}{\ln\left(k+\sqrt{k^2-1}\right) + \ln 2 + 1} > 0$, and proving irrationality for all these PCFs' limits. Note that these PCFs satisfy the convergence condition $A^2 + 4B > 0$ from Appendix A.
3. Additional values in Table 1.

These examples emphasize the strength of our approach: the determination of an irrationality measure $\delta$ without a need to find an explicit formula for a PCF's sequence or to find its limit.

## 2.3 Conditions for the existence of factorial reduction

Theorems 1 and 2 raise two critical questions: (1) which PCFs have FR and if so, then (2) what their exponential orders $\lambda$ are. Following many computer tests, the following conjecture is an effort to answer the first question. The second will be discussed later.

**Conjecture 1.1.** *Given a polynomial $b_n$, there exists an $a_n$ for which* $\mathrm{PCF}[a_n, b_n]$ *has FR if and only if $b_n$ is decomposable to linear factors over* $\mathbb{Q}$.

<u>Examples: The roots of $b_n$ and their effect on FR</u>

For the case where $b_n = n^2 - 2$, which possesses irrational roots, we did not discovere any $a_n$ for which $\mathrm{PCF}[a_n, b_n]$ has FR. Similarly, in scenarios where $b_n$ has nonreal roots, such as $b_n = n^2 + 1$, we observed the same outcome (see Appendix E).

In the case of $b_n = n^2 - 2$, which has irrational roots, we found no $a_n$ such that $\mathrm{PCF}[a_n, b_n]$ has FR, and similarly in the case where $b_n$ has nonreal roots, such as $b_n = n^2 + 1$ (see Appendix E).

On the other hand, we found numerically that for $b_n = 8n^2 - 2 = 2(2n + 1)(2n - 1)$, there exist choices of $a_n$ that provide PCFs with FR, for example:

$$\mathrm{PCF}[7n + 3, 8n^2 - 2] = 3 + \cfrac{6}{10 + \cfrac{30}{17 + \cfrac{70}{24 \ldots}}}. \qquad (\mathrm{Eq.}\ 4)$$

Another example for $b_n = -n^4 = -(n^2)(n^2)$ is

$$\mathrm{PCF}[2n^2 + 2n + 13, -n^4] = 13 - \cfrac{1}{17 - \cfrac{16}{25 - \cfrac{81}{37 \ldots}}}. \qquad (\mathrm{Eq.}\ 5)$$

Families of $a_n$s for which the PCFs have factorial reduction

The above $a_n$ choices are part of infinite families (see the following examples). In fact, our computer tests always find the $a_n$s to belong to infinite families that all have FR, and we propose the following conjectures.

**Conjecture 1.2.** *For each $b_n$ that is decomposable to linear factors over $\mathbb{Q}$, there exists at least one infinite family of $a_n$s for which every* $\mathrm{PCF}[a_n, b_n]$ *has FR.*

*Such a family is presented in* **Conjecture 1.3**.

This conjecture means that each $\mathrm{PCF}[a_n, b_n]$ with FR could be generalized to an infinite family of PCFs with the same $b_n$ and different $a_n$s. An interesting question is whether there would always exist members of this family that prove the irrationality of the constants to which they converge (see **Theorem 3** for more information).

Example of an infinite family of PCFs for $b_n = -n^4$ (found empirically)

The PCF in (Eq. 5) can be generalized, and considered as a special case of the following family of PCFs: For $b_n = -n^4$, and any $k$ of the form $k = m^2 - m + 1$, $m \in \mathbb{Z}$, the following PCF has FR:

$$\mathrm{PCF}[2n^2 + 2n + k, -n^4] = k - \cfrac{1}{4 + k - \cfrac{16}{12 + k - \cfrac{81}{24 + k \ldots}}}.$$

Example of an infinite family of PCFs for any general $b_n$ (found empirically)

**Conjecture 1.3.** *For any $d \in \mathbb{N}$, and $B, x_1, \ldots, x_{2d} \in \mathbb{Q}$, the following* $\mathrm{PCF}[a_n, b_n]$ *has FR:*

$$b_n = B \cdot \prod_{i=1}^{2d} (n - x_i)$$

$$a_n = \frac{B}{m} \cdot \prod_{i=1}^{d} (n + 1 - x_i) - m \cdot \prod_{i=d+1}^{2d} (n - x_i)$$

*for any* $m \in \mathbb{Q} \backslash \{0, \pm\sqrt{|B|}\}$ *(since then* $B = -A^2/4$, *or* $A = 0$*).*

The complete structure for $\deg b_n = 2$ and $\deg a_n = 1$

Having performed many numerical tests, we propose the next general conjecture for the $a_n$ families for which the PCF has FR for a given reducible $b_n$ of degree 2.

**Conjecture 1.4.** *For every* $b_n$ *of the form* $b_n = B(n - x_1)(n - x_2)$ *with* $B, x_1, x_2 \in \mathbb{Q}$, PCF$[a_n, b_n]$ *has FR if and only if* $a_n$ *belongs to one (or more) of the families*

$$a_n^{(k)} = \left(\frac{B}{m} - m\right) n + k, \qquad \text{(Eq. 6)}$$

$$a_n^{(k)} = k(2n - x_1 - x_2 + 1), \qquad \text{(Eq. 7)}$$

*for any* $k \in \mathbb{Q} \backslash \{0\}$, *and* $m \in \mathbb{Q} \backslash \{0, \pm\sqrt{|B|}\}$.

We emphasize that the above conjecture is formulated with rational parameters ($B, k$, and $m$), yielding rational PCFs. Alternatively, an equivalent conjecture can be formulated using integer parameters and yielding integer PCFs. The equivalence is achieved by multiplying $a_n$ and $b_n$ by a constant (see inflation process in Appendix C). For the sake of simplicity and coherence, the examples below are chosen to be integers ($B \in \mathbb{Z}$, $k \in \mathbb{Z}$ or $k \in \frac{1}{2}\mathbb{Z}$, and $m|B$). However, the conjecture is presented in the most general form, using rational numbers.

Numerical experiments show that these two families of PCFs, and only them, have FR for a given reducible $b_n$ of degree 2. Based on these numerical tests, we found exponentially tight ($\doteq$) formulas for $\mathrm{GCD}_n$ of some cases (Appendix E). Interestingly, the families share additional properties, e.g., their $\mathrm{GCD}_n$s are always found to be closely related (we do not yet understand the general structure).

<u>Example of the complete structure for $b_n = 8n^2 - 2$</u>

For $b_n = 8n^2 - 2$, the first family (Eq. 6) has the following PCFs for $m = 1$ and any $k \in \mathbb{Z}$ (note that (Eq. 4) is a special case with $k = 3$)

$$\mathrm{PCF}[7n + k, 8n^2 - 2] = 0 + k + \cfrac{6}{7 + k + \cfrac{30}{14 + k + \cfrac{70}{21 + k \ldots}}},$$

and for $m = 2$ and any $k \in \mathbb{Z}$

$$\mathrm{PCF}[2n + k, 8n^2 - 2] = 0 + k + \cfrac{6}{2 + k + \cfrac{30}{4 + k + \cfrac{70}{6 + k \ldots}}}.$$

The second family (Eq. 7) contains for any $k \in \mathbb{Z}$

$$\mathrm{PCF}[k(2n + 1), 8n^2 - 2] = k + \cfrac{6}{3k + \cfrac{30}{5k + \cfrac{70}{7k \ldots}}}.$$

We tested many of these PCFs numerically and indeed they all have FR, while other PCFs with the same $b_n = 8n^2 - 2$ and different $a_n$s do not.

## 2.4 Summary of our main conjectures regarding PCFs and factorial reduction

The above examples summarize the four aspects of our conjectures so far:

1. FR of PCF[ $\cdot$ , $b_n$]: an $a_n$ exists if and only if $b_n$ has rational roots only.
2. $a_n$ families: for every such $b_n$, there exist several infinite families of $a_n$s for which the PCFs have FR.
3. Each $a_n$ belongs to (at least) one of such families.
4. PCFs of the same family have closely related $\mathrm{GCD}_n$s.

A hint for the structure of arbitrary degrees

We expect similar results for any $b_n$ of degree $\geq 2$ and explain why. For any $b_n$, of any degree, we always find $a_n$s with leading coefficients $\left(\frac{B}{m} - m\right)$ for some $m$, same as the first family (Eq. 6) for $b_n$ of degree $= 2$. For example, PCF[$11n^2 - 2n - 1, 12n^4$], PCF[$4n^2 - 4n - 2, 12n^4$], and PCF[$n^2 - 6n - 3, 12n^4$] have FR, corresponding to $B = 12$ and $m = 1,2,3$, respectively. See Appendix E for more examples. We do not yet know what conditions must be satisfied by the other coefficients of the $a_n$s to have FR.

There exist other $a_n$s of a different and still unknown form, such as Apéry's and the example of **Conjecture 1.2** shown above. Furthermore, the form of the leading coefficient $A = \left(\frac{B}{m} - m\right)$ can be obtained by the following criterion: $A$ is of this form if and only if the equation $m^2 = A \cdot m + B$ has a rational root. It is interesting to try generalizing the above conjectures to discover the most general rules of this mathematical structure. Additional hints of the mathematical complexity of the yet unknown general structure are related to the existence of generalized Pythagorean triples (see Section 2.7 below).

### 2.5 Closed-form formula of $\mathrm{GCD}_n$ and the exponential order $\lambda$

The goal of the next sections is to predict exponentially tight formulas for $\mathrm{GCD}_n$s, i.e., up to a slower than exponential factor. For each $\mathrm{PCF}[a_n, b_n]$, we aim to find both the exponential order $\lambda$ and the closed-form expression for $\mathrm{GCD}_n$ that yields this $\lambda$. Representative examples are provided in Table 1 below. This table includes PCF examples of different sorts, some for which we found (numerically) exponentially tight formulas and others for which we did not.

| PCF | $0+\cfrac{2}{1+\cfrac{5}{2+\cfrac{10}{3\ldots}}}$ | $-15+\cfrac{250}{5+\cfrac{750}{25+\cfrac{1500}{45\ldots}}}$ | $2+\cfrac{2}{5+\cfrac{12}{8+\cfrac{30}{11\ldots}}}$ | $5+\cfrac{18}{14+\cfrac{55}{23+\cfrac{112}{32\ldots}}}$ | $1+\cfrac{4}{4+\cfrac{28}{7+\cfrac{70}{10\ldots}}}$ |
|---|---|---|---|---|---|
| $b_n$ | $n^2+1$ | $125n^2+125n$ | $4n^2-2n$ | $10n^2+7n+1$ | $9n^2-3n-2$ |
| $a_n$ | n | $20n-15$ | $3n+2$ | $9n+5$ | $3n+1$ |
| $\lambda$ | 0 - no FR | $0.5370 \approx \frac{5}{2\cdot\sqrt{3}\cdot e}$ | $0.5413 \approx 2\cdot\frac{2}{e^2}$ | $0.2180 \approx ?$ | 3 |
| $\mathrm{GCD}_n \doteq$ | $\mathrm{GCD}_n$'s growth rate<br>is less than a factorial | $n!\cdot\frac{5^n}{2^n 3^{\lfloor n/2\rfloor}\mathrm{LCM}[n]}$ | $(2n-1)!!\cdot\frac{2^n}{\mathrm{LCM}[2n]}$<br>!! is double factorial | ? | $(3n+1)!!!$<br>!!! is triple factorial |
| $\delta$ | $-1$ | -0.58 | -0.31 | -0.4 | 1 |
| Remarks | $b_n$ has unreal roots | $5^n$ is due to<br>inflation;<br>see Appendix C | Equivalent<br>representation:<br>$\mathrm{GCD}_n \doteq$<br>$n!\cdot\frac{\binom{2n}{n}}{\mathrm{LCM}[2n]}$ | "Zebra" pattern<br>(see below).<br>2, 5, and 11 are<br>exponentially<br>coprime to $\mathrm{GCD}_n$ | This PCF is an<br>inflation of the<br>*regular* continued<br>fraction of $\varphi$;<br>see Appendix C |

**Table 1: PCF examples with different $\mathrm{GCD}_n$ formulas.** The presented $\mathrm{GCD}_n$ formulas differ from the exact $\mathrm{GCD}_n$ by sub-exponential factors.

The term "exponentially coprime to $\mathrm{GCD}_n$" generalizes the idea of a coprime and means that the highest powers of $p$ dividing $\mathrm{GCD}_n$ for the terms in the sequence increase sub-exponentially. This statement implies that a certain prime does not affect the reduction.

Generalizing from these examples and many more (Appendix E), a conjectured structure of the exact $\mathrm{GCD}_n$ forms (up to sub-exponential factors) is presented next. Note that there exist multiple equivalent ways to present some of the forms, for example using $(2n-1)!! \cdot \frac{2^n}{\mathrm{LCM}[2n]} = n! \cdot \frac{\binom{2n}{n}}{\mathrm{LCM}[2n]}$ (as shown in Table 1). Also, the notation $f_n = \Theta(g_n)$ mean that $f$ is bounded both above and below by $g$ asymptotically, i.e., there exist positive constants $c, C \in \mathbb{R}^+$ such that $c \cdot g_n \leq f_n \leq C \cdot g_n$ for large enough $n$s.

**Conjecture 2** *(The exact forms of* $\mathrm{GCD}_n$*). We can represent every* $\mathrm{GCD}_n$ *as a multiple of two parts, a factorial and an exponential expression:*

- *The factorial part in general appears in the form* $\prod_{i=1}^{n} P(i)$*, where* $P$ *is an integer polynomial of degree* $\deg a_n$

  (*The special case of* $n!$ *corresponds to* $P(x) = x$ *and* $\deg(a_n) = 1$).
- *The exponential part takes one of the following forms or their multiples:*
  - *Power of a prime:* $p^{\Theta(n)}$
    - *In the numerator: We found only integer or half-integer powers, such as* $5^n$ *and* $7^{\lfloor n/2 \rfloor}$*. In some PCFs, these powers can be explained by inflation (see Appendix C).*
    - *In the denominator: Only primes* $p$ *raised to the power of* $\left\lfloor \frac{n}{p-1} \right\rfloor \deg(a_n)$*, such as* $11^{\left\lfloor \frac{n}{10} \right\rfloor}$ *and* $2^n$*. Note that this exponent* $p^{\left\lfloor \frac{n}{p-1} \right\rfloor}$ *conforms to the highest*

*power of* $p$ *that divides the factorial expression* $(un+v)!^{(u)}$ *(when* $u$ *and* $p$ *are coprime). Hence, when this exponent appears in the denominator,* $\mathrm{GCD}_n$ *is exponentially coprime with* $p$*.*

- *When part of the expression is a "Zebra" (see below), we find more complicated fractional powers in the numerator, such as* $2^{\lfloor 3n/4 \rfloor}$*.*

- $\mathrm{LCM}[f \cdot n]$ *for some* $f \in \mathbb{N}$*: For example, the Apéry's work has* $\mathrm{LCM}[n]$*, and Table 1 shows a case of* $\mathrm{LCM}[2n]$*. This is seen only in denominators.*
- *Zebra: There is an additional pattern for which we lack an explicit formula. We find this pattern in the denominators. We can identify this pattern in many PCFs but do not entirely understand it. The investigation of the Zebra pattern is left to future work.*

For computational simplicity, most of our numerical analysis is focused on PCFs of $\deg b_n = 2$ and $\deg a_n = 1$. Based on this analysis and additional simulations, we conjecture that the above description captures any $\mathrm{GCD}_n$ sequence of a PCF, also of the higher order $a_n, b_n$. Furthermore, we expect analogous mathematical structures to exist in greatest common divisors of any linear recursion with polynomial coefficients, the investigation of which remains for future work. Note that, for PCFs without FR, numerical analyses show that $\Theta\big(\pi(n)\big)$ primes are exponentially coprime to $\mathrm{GCD}_n$, where $\pi(n)$ is the prime-counting function.

### 2.6 Fast calculation of PCFs using simplified recursion formulas and factorial reduction

In this section, we discuss an application of the ability to predict the exact formulas of FR and other forms of reduction. Provided we have a closed-form formula for $\mathrm{GCD}_n$, we can apply the reduction straight to the recursion, so that the computation is performed with smaller integer values. Such simplified recursions enable faster estimation of the PCF limit. The computation advantage of such recursion is substantial with FR: it requires only manipulating sequences that grow exponentially with the PCF depth (instead of super-exponentially).

Example: A simple recursion for the reduced numerator and denominator

For $b_n = 2n^2 + n, a_n = n$, we (numerically) find that $\mathrm{GCD}_n \doteq \frac{n!}{2^n}$, and therefore, there exist integer sequences $p_n'$ and $q_n{}'$ such that

$$p_n = \frac{n!}{2^n} \cdot p_n'$$

$$q_n = \frac{n!}{2^n} \cdot q_n'.$$

Thus, $\mathrm{GCD}[p_n', q_n'] \doteq 1$ so that the majority of the original $\mathrm{GCD}_n$ is being used in this reduction. We can substitute in the recursion (Eq. 1) that $p_n$ and $q_n$ both uphold, and yield a recursion for $p_n'$ and $q_n'$ (after simple manipulations)

$$n(n-1)u_n' \ = 2n\, u_{n-1}' + 4u_{n-2}',$$

which has rational coefficients. This recursion generates the reduced numerator and denominator sequences. In fact, for **any** integer initial values, this recursion generates an **integer** sequence (or rational with sub-exponential denominators). Additional pairs of polynomials $a_n, b_n$ with the same properties are presented in Appendix E.

### 2.7 Hints for a deeper mathematical structure

This section provides additional examples of special mathematical properties that we found numerically and hint at a much wider theory that still awaits discovery.

For every $b_n^{(x)} = x^2(2n^2 + n), x \in \mathbb{Z}$, there exist families of $a_n$ for which PCF$[a_n, b_n]$ has FR. These $a_n$ families include

$$a_n = z^{(x)} \cdot n + k,$$

where $k$ is an integer, but the options for integer $z^{(x)}$ are finite. These $z^{(x)}$s are precisely the integers for which there exists a **(-6)-Pythagorean triple** $(x, y, z)$ for some $y$; i.e., the Diophantine equation $z^2 = x^2 + y^2 + 6xy$ has a solution. This was discovered with the help of OEIS [27]. This result is a special case that coincides with the general structure we discovered for $\deg b_n = 2$ (see **Conjecture 1.4**).

Apéry wrote two more pairs of polynomials, the PCFs of which prove the irrationality of $\zeta(2)$ and $\ln 2$. After considering **Conjecture 1.2** (infinite families for a given $b_n$), we looked for these families with the others $a_n$s. For $\ln 2$, where $b_n = -n^2$, we discovered $a_n = k(2n + 1), k \in \mathbb{Z}$ as the particular structure for $\deg b_n = 2$ predicts. Moreover, for odd $k$s, such as the Apéry's ($k = 3$), we get $\mathrm{GCD}_n \doteq \frac{n!}{\mathrm{LCM}[n]}$, and for even $k$s, $\mathrm{GCD}_n \doteq \frac{n!}{2^n \mathrm{LCM}[n]}$. A generalized proof for this case, even without knowing the PCFs limits, is available in Appendix D. The theorem in the next section shows how almost any $k$ constructs a PCF that proves the irrationality of its limit, although, apart from $\ln 2$, the identity of these irrational limits is still unknown to us.

### 2.8 Infinite $a_n$s that prove irrationality for a given $b_n$

This next section shows infinite families of PCFs that prove irrationality of certain numbers. Specifically, we conjecture that for any $b_n$, there exists an infinite set of $a_n$s such that each constructs a PCF that proves the irrationality of its limit.

**Theorem 3.** *We consider families of balanced-degree* PCFs *that have FR and are created from a constant* $b_n$*, and* $a_n^{(k)}$ *that are multiples of a single polynomial, i.e.,* $a_n^{(k)} \in \left\{k \cdot a_n^{(1)} \mid k \in \mathbb{Z}\right\}$. *Assuming that* $\lambda$*, the exponential orders of* $\mathrm{GCD}_n$*, is bounded as a sequence in* $k$ *(based on part 4 of* ***Conjecture 1****'s summary), we find*

$$\lim_{k\to\infty} \delta = 1.$$

*In particular, for large enough* $k$*s, the limit will be irrational since* $\delta > 0$.

*Proof.* (Straightforward) If $k \to \infty$, then the leading coefficient of $a_n^{(k)}$ uphold $A_k \to \infty$ and the characteristic equation $\alpha^2 = A_k \cdot \alpha + B$ has a solution that certifies $\alpha_k \to \infty$. Substituting in **Theorem 2** with a constant $B$, while assuming $\ln \lambda$ is bounded, we have

$$\delta = \frac{\ln|\alpha| - \ln|B| + \ln \lambda}{\ln|\alpha| - \ln \lambda} \to \frac{\ln|\alpha|}{\ln|\alpha|} = 1. \qquad \square$$

<u>Observation</u>: Combining this theorem and **Conjecture 1.4**, we expect that, for any $b_n$ of degree 2 with rational-only roots, there exists an infinite set of $a_n$s such that $\mathrm{PCF}[a_n, b_n]$ proves the irrationality of its limit. As for higher degrees, we conjecture the existence of similar structures.

## 2.9 Additional properties of the greatest common divisor

We investigate additional results that can help prove properties about $\mathrm{GCD}_n$, for all PCFs cases, either with or without FR. Thus far in our paper we analyzed the growth rate of the sequence $\mathrm{GCD}_n$ as a function of the PCF's depth, $n$. One property that we examined and believe could be useful for proving some of our conjectures is whether $\mathrm{GCD}_n$ divides its consecutive $\mathrm{GCD}_{n+1}$. We find numerically that in many of the fractions this does not hold for the definition of $\mathrm{GCD}_n$ sequence, i.e., $\mathrm{GCD}_n \nmid \mathrm{GCD}_{n+1}$ for some $n$s.

The above observation motivates the study of the greatest common divisor of two consecutive numerators and denominators, defined as

$$\mathrm{GCD2}_n \stackrel{\mathrm{def}}{=} \mathrm{GCD}[p_n, q_n, p_{n-1}, q_{n-1}].$$

By this definition and the recursion formula $(\mathrm{Eq.}\,1)$ for $p_n$ and $q_n$, one can show that for all $n$

$$\mathrm{GCD2}_n \mid \mathrm{GCD2}_{n+1}.$$

Since $\mathrm{GCD2}_n \mid \mathrm{GCD}_n$, part of the reduction may be explained by $\mathrm{GCD2}_n$. It remains to be seen what part of the FR and its exponential part is contained in $\mathrm{GCD2}_n$. Having inspected many PCFs numerically, with or without FR, we conjecture the following.

**Conjecture 3.** *For any* PCF*:*

$$\mathrm{GCD2}_n \doteq \mathrm{GCD}_n,$$

The fact that both $\mathrm{GCD}_n$ and $\mathrm{GCD2}_\mathrm{n}$ exhibit the same exponential order means that all the theorems and conjectures presented here may also apply to $\mathrm{GCD2}_n$. Specifically, if FR exists, then both the factorial and the exponential part of $\mathrm{GCD}_n$ will exist in $\mathrm{GCD2}_n$. The important consequence is that we can use either $\mathrm{GCD}_n$ or $\mathrm{GCD2}_n$ for purposes of irrationality proofs, as **Theorem 2**.

This conjecture enables us to treat $\mathrm{GCD2}_n$ as a growing product of some integer series and, at a given depth $n$, calculate and reduce only one integer term: $\frac{\mathrm{GCD2}_n}{\mathrm{GCD2}_{n-1}}$. For example, if $\mathrm{GCD2}_n = \frac{n!}{\mathrm{LCM}[n]}$, we can reduce the numerators and the denominators **at each depth** $n$ by the factor $n/p$ if $n$ is a power of some prime number $p$ and by $n$ if it is not.

Moreover, this definition is advantageous because $\mathrm{GCD2}_n|\mathrm{GCD}_n$ for all $n$, and it thus sorts out sub-exponential factors that have no effect on proving irrationality. This observation facilitates the numerical analysis and helps identify the exact formula for the $\mathrm{GCD2}_n$.

As a side note, **Conjecture 3** helps show that $\mathrm{GCD}_n$s of PCFs that have FR always have a factorial term such as $(n!)^d$ for some integer $d$, rather than a term such as $n^{d \cdot n}$ (which also grows like $(n!)^d$ up to an exponential factor by Stirling's approximation). In fact, all the PCFs with FR that we encountered could be written as $(n!)^d \cdot \frac{S_n}{R_n}$ with $S_n$ and $R_n$ being integer sequences that grow exponentially. Some cases are more complex, such as when $\mathrm{GCD}_n \doteq (3n+1)!!!$, but these do not contradict the above statement. It would be interesting to try to prove this phenomenon.

# 3 Discussion and Open Questions

## 3.1 Outlook and motivation

By their further development, the conjectures presented can provide useful tools for irrationality proofs, as well as for fast calculations of polynomial integer recursions of mathematical constants.

Specifically, the results related to FR can be applied to shrink the search space of the Ramanujan Machine algorithms [10,30]. By focusing on PCFs with FR, the algorithms would have a better chance to find new conjectures that are simultaneously of a relatively fast computation time and have nontrivial $\delta$s that we can extract. That is, removal of the cases that have no FR avoids all the hard-to-compute PCFs that also provide trivial $\delta$s.

Looking forward, we believe that by generalizing the mathematical structure of PCFs with FR, it would be possible to find universal structures in PCFs made from arbitrary-degree polynomials. As a more ambitious step, it is interesting to consider deeper linear recursions (beyond depth 2), which can also be harnessed to find new conjectures. One can search for analogous algebraic structures and ideas as presented above.

In the following, we present several ideas and open questions that arise from our mathematical experiments and from our conjectures. These open questions may be simple or complex, and we hope that they can engender more ideas for future research in different communities.

### 3.4 Implications of factorial reduction for a faster computation of PCFs

Once a closed formula for $\mathrm{GCD}_n$ has been found, numerical calculations of PCFs will become easier and faster since the FR considerably decreases the numbers participating in the arithmetic operations. In particular, PCFs with FR benefit greatly from this reduction since $p_n$ and $q_n$ decrease from a super-exponentially (factorial) growth to exponential growth. In other words, finding the exact formula for the reduction enables one to construct a simpler recursion formula that directly gives the reduced numerators and denominators.

### 3.5 Families of PCFs

Following **Conjecture 1.2**, it is natural to try to generalize the families of $a_n$,$b_n$ for higher degrees. What affects the number of families and subfamilies? **Conjecture 1.4** claims that only two families exist for the discussed degrees, and one of them is branched into several subfamilies. In this case, the number of subfamilies depends merely on the number of divisors of $B$ (the leading coefficient of $b_n$). We do not yet have a solid and more general conjecture that relates to all degrees. Another question regarding families of $a_n$ or $b_n$ is whether a relation exists between the limits of any sibling PCFs. For example, if this relation hints that the limits are equivalent, for proofs of irrationality, it will suffice to find just one limit and use **Theorem 3** (infinite $a_n$s that proves irrationality).

### 3.6 Finding and proving the exact form of $\mathbf{GCD}_n$

We did not find the exact form of the $\mathrm{GCD}_n$, but nevertheless tried to list the different types of expressions that comprise it. The motivation to find the closed form of $\mathrm{GCD}_n$ is the possibility of writing a reduced recursion that yields the reduced numerators and denominators, which can simplify any numerical calculation of the PCF. Moreover, a closed-form formula would also directly predict the bound $\delta$ given by the PCF.

We note that Apéry proved his case by finding an explicit expression for the PCF at each depth $n$, i.e., $p_n/q_n$. As an example of taking a more general approach, in Appendix D we address various $\mathrm{GCD}_n$s and bypass the need for an explicit expression. As examples that can promote future research, we present in Appendix E a set of unproven examples that precisely yield same simplified recurrence relations.

### 3.5 Predicting the exponential order $\lambda$:

To search for conjectures in the form of PCFs that prove the irrationality of constants, it suffices to predict only the exponential order $\lambda$. Using this value, **Theorem 2** calculates the bound $\delta$ on the irrationality measure. It remains to find a direct relation from $a_n$, $b_n$ to $\lambda$.

Acknowledgment

Special thanks should be given to Dr. Yaar Solomon and Prof. Yair Glasner from the Ben-Gurion University of the Negev for their mathematical guidance and constructive recommendations on writing this paper.

conflict of interest

The authors declare no conflict of interest.

# 5 Appendix

## Appendix A: Classification of PCFs

All PCFs can be split into three cases by the ratio of the degrees of the polynomials $a_n, b_n$ (Table 2). This ratio determines the PCF's convergence rate and the growth rate of $p_n$ and $q_n$. In most of this paper, we focused on Case 2, balanced-degree PCFs, where $\frac{d_b}{d_a} = 2$ (3rd column in the table below). PCFs of this case are those that prove the irrationality of $\zeta(3)$, and $\zeta(2)$ in the Apéry's work [1,2], as well as many more mathematical constants.

In this section, the other cases are presented and analyzed to ensure completeness and to lay the groundwork for future research. Interestingly, the analysis of the three cases demonstrates that not only Case 2 but also Cases 1 and 3 can provide expressions for $\delta$ used for irrationality proofs. These cases are generally less interesting because $\delta$ is independent of the size of $\mathrm{GCD}_n$, and thus provide proofs of irrationality only under stricter conditions. This comparison underscores the importance of PCFs of Case 2, where the size of $\mathrm{GCD}_n$ (both factorial and exponential components) significantly influences $\delta$ and its analysis is pivotal in the pursuit of irrationality proofs. Refer to the 5th row in Table 2 for further details.

| | Case 1 | Case 2 – balanced-degree | Case 3 |
|---|---|---|---|
| Degree ratio | $\frac{\deg b_n}{\deg a_n} > 2$ | $\frac{\deg b_n}{\deg a_n} = 2$ | $\frac{\deg b_n}{\deg a_n} < 2$ |
| Convergence[1] (see Appendix B, and Appendix D of [30]) | Conjectured if $B > 0$, and some cases[1] of $B < 0$ | If $A^2 + 4B > 0$, and some cases[1] of equality | Always |
| Rate of convergence (see Appendix B) | Conjectured to be Polynomial | Exponential | Super-exponential |
| Examples (from [10,30]) | $\frac{4+\pi}{\pi} = 2 + \cfrac{1}{2 + \cfrac{9}{2 + \cfrac{25}{2 \ldots + \cfrac{(2n-1)^2}{2 \ldots}}}}$ | $\frac{4}{\pi} = 1 + \cfrac{1}{3 + \cfrac{4}{5 + \cfrac{9}{7 \ldots + \cfrac{n^2}{2n+1 \ldots}}}}$ | $\frac{1+e}{-1+e} = 2 + \cfrac{1}{6 + \cfrac{1}{10 + \cfrac{1}{14 \ldots + \cfrac{1}{4n+2 \ldots}}}}$ |
| Does this PCF provide a **nontrivial** $\delta$? (see Appendix B) | Only if there is a reduction of $q_n$ after which $\frac{q_n}{\mathrm{GCD}_n}$ is **sub-**exponential.<br>Namely, it is necessary that<br>$\mathrm{GCD}_n \doteq q_n$<br>but it is not sufficient<br>(This is a conjecture) | If and only if there is FR (see **Theorem 1**).<br>$\delta$ size depends on the $\mathrm{GCD}_n$'s exponential size, following **Theorem 2** | Always provides<br>$\delta = 1 - \frac{\deg b_n}{\deg a_n}$<br>or better if there is FR.<br>Proves irrationality even without FR if<br>$\deg b_n < \deg a_n$[2] |

**Table 2: Summary of the three cases of PCFs**, partitioned by the ratio of the degrees of $a_n, b_n$. We show the conditions for each case that leads to a nontrivial bound $\delta$ on the irrationality measure of the limit. Note the crucial role of FR in the balanced-degree type (Case 2), which is at the core of this manuscript. $A$ and $B$ are the leading coefficients of $a_n$ and $b_n$, respectively. [1]There are cases for which we do not know the conditions for convergence. Also, note that PCFs might converge to infinity (see Lemma 4 in Appendix B). [2]This condition coincides with Tietze's criterion [26].

## Appendix B: Proof for Theorems 1 and 2

The goal of this appendix is to provide a provide a proof for Theorems 1 and 2 from the main text. For convenience, we collect our notation at this point.

Let $a, b \in \mathbb{Z}[X]$ be two fixed polynomials of degrees $d_a, d_b$ and leading coefficients $A, B$, respectively. We write $a_n = a(n)$ and $b_n = b(n)$ for any integer $n$, and consider the polynomial continued fraction $\mathrm{PCF}[a_n, b_n]$ as introduced in Section 1.2. The numbers $p_n$ and $q_n$ denote the numerator and denominator of the $n$-th convergent $p_n/q_n$ of $\mathrm{PCF}[a_n, b_n]$, satisfying the recursion $u_n = a_n u_{n-1} + b_n u_{n-2}$. We will assume that $b_n \neq 0$ for all $n \in \mathbb{N}$ to guarantee that the PCF is infinite.

We also denote by $L \overset{\text{def}}{=} \lim_{n\to\infty} \frac{p_n}{q_n}$ the limit of the PCF – assumed to exist for the sake of Theorems 1 and also finite for Theorem 2. For integer sequences $p_n, q_n$ generating a rational sequence $p_n/q_n \neq L$ converging to $L$, we define

$$\delta_{p_n,q_n} \overset{\text{def}}{=} \sup\left\{\delta \in [-1,\infty] \,\middle|\, \left|L - \frac{p_n}{q_n}\right| < \frac{1}{|q_n|^{1+\delta}} \text{ for infinitely many } n\right\}.$$

We will prove a lower bound on the irrationality measure of $L$ by bounding the above expression.

Without loss of generality, we will assume that $A > 0$. Otherwise, one could consider $\mathrm{PCF}[-a_n, b_n]$, whose limit is $-L$ and convergents are $-p_n/q_n$, as can be proven using inflation (see Appendix C) by the constant sequence $-1$.

**Step 1: A lower bound on the irrationality measure**

The first step towards proving our theorems consists of producing a controllable lower bound on $\delta_{p'_n,q'_n}$ where $(p'_n, q'_n) \stackrel{\text{def}}{=} \left(\frac{p_n}{\mathrm{GCD}_n}, \frac{q_n}{\mathrm{GCD}_n}\right)$.

**Lemma 1.** *With the notation above, we have the inequality*

$$\delta_{p'_n,q'_n} \geq \delta \stackrel{\text{def}}{=} \limsup_{n\to\infty} \frac{\ln|q_{n+1}| - \ln\left|\prod_{i=1}^{n} b_i\right| + \ln \mathrm{GCD}_n}{\ln|q_n| - \ln \mathrm{GCD}_n}.$$

This is a consequence of the subsequent lemma.

**Lemma 2.** *The following inequality holds for sufficiently large* $n$

$$\left|L - \frac{p'_n}{q'_n}\right| = \left|L - \frac{p_n}{q_n}\right| \leq \left|\frac{\prod_{i=1}^{n} b_i}{q_{n+1} q_n}\right|.$$

*Proof.* See [29] at Appendix D.1. □

*Proof of Lemma 1*. Observe that the expression defining $\delta$ can be obtained by

$$\delta = \sup\left\{ d \in [-1,\infty] \;\middle|\; \left|\frac{\prod_{i=1}^{n} b_i}{q_{n+1} q_n}\right| < \frac{1}{\left|\frac{q_n}{\mathrm{GCD}_n}\right|^{1+d}} \text{ for infinitely many } n \right\}.$$

$$= \sup\left\{ d \in [-1,\infty] \;\middle|\; d < \frac{\ln|q_{n+1}| - \ln\left|\prod_{i=1}^{n} b_i\right| + \ln \mathrm{GCD}_n}{\ln|q_n| - \ln \mathrm{GCD}_n} \text{ for infinitely many } n \right\}.$$

In combination with Lemma 2, this equality establishes $\delta_{p'_n,q'_n} \geq \delta$ . □

Also, we conjecture this lower bound to be tight, suggesting that the largest lower bound provided by this PCF on the irrationality measure of $L$ is the $\delta$ we defined.

**Step 2: Exponentially tight estimation for $\prod_{i=1}^{n} b_i$**

**Lemma 3.** *For any polynomial* $b_n$

$$\prod_{i=1}^{n} b_i \doteq B^n \cdot n!^{d_b}.$$

*Proof.* We claim that $\prod_{i=1}^{n} b_i$ and $B^n \cdot n!^{d_b}$ differ by factors that grow (or decay) more slowly than exponentials, and possibly a constant sign. To show that claim, assume first that $B > 0$, and define some positive integer $k$ such that

$$|b_n - B \cdot n^{d_b}| \leq k \cdot B \cdot n^{d_b - 1},$$

$$|a_n - A \cdot n^{d_a}| \leq k \cdot A \cdot n^{d_a - 1},$$

for all $n \in \mathbb{N}$. Then, for $n > k$, we have

$$0 < B \cdot n^{d_b} \left(\frac{n-k}{n}\right) \leq b_n \leq B \cdot n^{d_b} \left(\frac{n+k}{n}\right), \qquad \text{(Eq. B1)}$$

$$0 < A \cdot n^{d_a} \left(\frac{n-k}{n}\right) \leq a_n \leq A \cdot n^{d_a} \left(\frac{n+k}{n}\right) \qquad \text{(Eq. B2)}$$

(we will use (Eq. B2) later). Thus,

$$\frac{\prod_{i=1}^{n} b_i}{B^n \cdot n!^{d_b}} \leq \frac{\prod_{i=1}^{k} b_i}{B^k k!} \prod_{i=k+1}^{n} \frac{i+k}{i} = \frac{\prod_{i=1}^{k} b_i}{B^k k!} \frac{k!}{(2k)!} \cdot \frac{(n+k)!}{n!} = \Theta(n^k),$$

$$\frac{\prod_{i=1}^{n} b_i}{B^n \cdot n!^{d_b}} \geq \frac{\prod_{i=1}^{k} b_i}{B^k k!} \prod_{i=k+1}^{n} \frac{i-k}{i} = \frac{\prod_{i=1}^{k} b_i}{B^k k!} k! \cdot \frac{(n-k)!}{n!} = \Theta(n^{-k}).$$

These bounds show that $\prod_{i=1}^{n} b_i$ and $B^n \cdot n!^{d_b}$ differ by (up to) polynomial factors, and perhaps a constant sign, as required for the relation "$\doteq$", defined at (Eq. 3).

Note that minor adjustments are required for establishing similar bounds to (Eq. B1) when $B < 0$. □

**Step 3: Exponentially tight estimation for $q_n$**

Our first step towards estimating the behavior of $q_n$ as $n \to \infty$ will be to understand the asymptotic behavior of the quotient $q_n/q_{n-1}$. The situation for this quotient is particularly simple in the case of constant coefficients $a'_n = A'$ and $b'_n = B'$, which we briefly describe here.

**Example.** The recursion $q'_n = A'q'_{n-1} + B'q'_{n-2}$ can be re-written in a matrix form as

$$\begin{pmatrix} q'_n \\ q'_{n-1} \end{pmatrix} = \underbrace{\begin{pmatrix} A' & B' \\ 1 & 0 \end{pmatrix}}_{T} \begin{pmatrix} q'_{n-1} \\ q'_{n-2} \end{pmatrix} = T^n \begin{pmatrix} q'_0 \\ q'_{-1} \end{pmatrix}.$$

In terms of the Möbius action, this translates to the equality

$$\frac{q'_n}{q'_{n-1}} = T^n \left( \frac{q'_0}{q'_{-1}} \right) \qquad \text{(Eq. B3)}$$

for every $n$. Observe now that the characteristic polynomial of $T$ is $x^2 - A'x - B'$, and that

$$T = P \begin{pmatrix} \lambda_- & 0 \\ 0 & \lambda_+ \end{pmatrix} P^{-1} \quad \text{for} \quad P = \begin{pmatrix} \lambda_- & \lambda_+ \\ 1 & 1 \end{pmatrix},$$

where $\lambda_\pm = (A' \pm \sqrt{A'^2 + 4B'})/2$ are the roots of the characteristic polynomial. In virtue of the assumptions $A' > 0$ and $A'^2 + 4B' > 0$, the roots are real, distinct, and $|\lambda_+| > |\lambda_-|$. Combining this diagonalization statement with (Eq. B3), we conclude that

$$\lim_{n \to \infty} q'_n/q'_{n-1} = \lambda_+$$

for all initial conditions $(q_0', q_{-1}')$, except if $\frac{q_0'}{q_{-1}'} = \lambda_-$, in which case $\lim_{n\to\infty} q_n'/q_{n-1}' = \lambda_-$.

In words, the example above shows that the quotient $q_n'/q_{n-1}'$ converges to the dominant (larger in absolute value) eigenvalue $\lambda_+$ of the characteristic polynomial for generic initial conditions $(q_0', q_{-1}')$, under the assumption of constant coefficients. The following theorem by Poincaré [31] states that an similar behavior is also expected when the coefficients are not necessarily constant, yet admit finite limits as $n \to \infty$.

**Theorem A.** *(Poincare 1885) Let* $u_n = a_n' u_{n-1} + b_n' u_{n-2}$, $n \in \mathbb{N}$, *be a linear recurrence with real-valued coefficients* $a_n', b_n'$ *and initial conditions* $(u_0, u_{-1}) \in \mathbb{R}^2 \setminus \{(0,0)\}$. *Assume that* $a_n', b_n'$ *have limits* $\lim_{n\to\infty} a_n' = A'$, $\lim_{n\to\infty} b_n' = B'$ *with* ${A'}^2 + 4B' > 0$, *and* $\lambda_+$, $\lambda_-$ *are two distinct real solutions of the equation* $x^2 - A'x - B' = 0$. *Then, every non-zero solution* $u_n$ *of the recurrence fulfills*

$$\lim_{n\to\infty} \frac{u_n}{u_{n-1}} = z \in \{\lambda_+, \lambda_-\}. \qquad \text{(Eq. B4)}$$

*Moreover, if* $|\lambda_+| > |\lambda_-|$, *there exists at most one number* $c \in \mathbb{R} \cup \{\infty\}$ *such that*

$$z = \lambda_+$$

*for all initial conditions* $(u_0, u_{-1})$ *except if* $u_0/u_{-1} = c$.

*About the proof.* (Eq. B4) is a theorem due to Poincaré [31] (see also [32] section 2.14 as a secondary reference). The remaining statement follows from (Eq. B4), the equality

$$\frac{u_n}{u_{n-1}} = S_n \left(\frac{u_0}{u_{-1}}\right) \quad \text{with} \quad S_n = T_n \cdots T_1 \ \text{ and } \ T_k = \begin{pmatrix} a_k' & b_k' \\ 1 & 0 \end{pmatrix}$$

for every $n$, and a theorem by Piranian and Thron [32, Theorem. 1], regarding the pointwise limit of sequences of Möbius transformations, applied to the sequence $S_n$. □

Relying on Theorem A, we will now prove the lemma below. Following the reasoning at Appendix A, we split in its statement all PCFs into three cases by the ratio $\frac{d_b}{d_a}$ of the degrees of the polynomials.

**Lemma 4.** *Let $\alpha$ be the solution of the equation $x^2 - Ax - B = 0$ with larger absolute value, assume that $A^2 + 4B > 0$ if $d_b/d_a = 2$, and that $B > 0$ if $d_b/d_a > 2$. Then*

$$q_n \doteq \begin{cases} B^{\frac{n}{2}} \cdot n!^{\frac{d_b}{2}} & \text{if} \quad d_b/d_a > 2 \quad \text{(Case 1)} \\ \alpha^n \cdot n!^{d_a} & \text{if} \quad d_b/d_a = 2 \quad \text{(Case 2)} \\ A^n \cdot n!^{d_a} & \text{if} \quad d_b/d_a < 2 \quad \text{(Case 3)} \end{cases}\cdot$$

*Proof.* Let $d \stackrel{\text{def}}{=} \max\left\{d_a, \frac{d_b}{2}\right\}$. Applying the substitution $q_n' = \frac{q_n}{n!^d}$ to the original recurrence, we obtain the new recurrence

$$q_n' = a_n'\, q_{n-1}' + b_n' q_{n-2}'$$

where $a_n' \stackrel{\text{def}}{=} \frac{a_n}{n^d}$ and $b_n' \stackrel{\text{def}}{=} \frac{b_n}{n^d (n-1)^d}$ for $n \geq 2$ (see inflation process at Appendix C), and $q_{-1}, q_0, q_1$ are as in the original recurrence. Observe that the new recurrence satisfies the condition of Theorem A, and the dominant eigenvalue $\lambda_+$ equals

$$\lambda_+ = \begin{cases} B^{1/2} & \text{if} \quad d_b/d_a > 2 \quad \text{(Case 1)} \\ \alpha & \text{if} \quad d_b/d_a = 2 \quad \text{(Case 2)}. \\ A & \text{if} \quad d_b/d_a < 2 \quad \text{(Case 3)} \end{cases}$$

Then, (Eq. B4) yields

$$\frac{q_n'}{q_{n-1}'} \to z \quad \text{as} \quad n \to \infty.$$

This implies that for every $\epsilon > 0$, there exist positive constants $r, R \in \mathbb{R}^+$ and $\kappa \in \{-1,1\}$ such that for large enough $n$

$$r(1-\epsilon)^n \leq \kappa \frac{q_n'}{z^n} = \kappa \frac{q_n}{z^n \cdot n!^d} \leq R(1+\epsilon)^n.$$

If $z = \lambda_+$, then the above proves that $\ln(1-\epsilon) \leq \lim_{n\to\infty} \frac{1}{n} \ln \kappa \frac{q_n}{\lambda_+^n \cdot n!^d} \leq \ln(1+\epsilon)$ for any $\epsilon > 0$, implying $q_n \doteq \lambda_+^n$ and establishing the lemma. On the other hand, if $z = \lambda_-$, then, given $c = q_0/q_{-1} \neq p_0/p_{-1}$, we find $\frac{p_n}{q_n} \doteq \left(\frac{\lambda_+}{\lambda_-}\right)^n$, resulting in a PCF that does not have a finite limit, which was excluded from our discussion. □

**Step 4: Combining these results with the irrationality criterion**

Now, for each case, we insert the estimation we found for $q_n$ into the expression

$$\delta = \limsup_{n\to\infty} \frac{\ln|q_{n+1}| - \ln\left|\prod_{i=1}^n b_i\right| + \ln \mathrm{GCD}_n}{\ln|q_n| - \ln \mathrm{GCD}_n}.$$

For convenient, define for some sub-exponential pre-factor as $E_n^b \stackrel{\text{def}}{=} \frac{\prod_{i=1}^n b_i}{B^n n!^{d_b}} \doteq 1$, allowing us to use the exact equality $\prod_{i=1}^n b_i = E_n^b \cdot B^n n!^{d_b}$ (and not just "$\doteq$"). Similarly, define $E_n^{\mathrm{GCD}}$ for $\mathrm{GCD}_n$ and $E_n^q$ for $q_n$ with its exponentially tight estimations at each case from **Lemma 4**. Also, note that $\ln n! \in \Theta(n \ln n)$.

Case 1: $\frac{d_b}{d_a} > 2$, $q_n \doteq B^{\frac{n}{2}} \cdot n!^{\frac{d_b}{2}}$

Note that

$$\left|\frac{\prod_{i=1}^n b_i}{q_{n+1} q_n}\right| \doteq \frac{B^n \cdot n!^{d_b}}{\left(B^{\frac{n}{2}} \cdot n!^{\frac{d_b}{2}}\right)^2} \doteq 1,$$

so we do not have a bound on the convergence rate. By numerical tests, we conjecture the converge rate to be polynomial. For this reason, to provide a nontrivial $\delta$, $\mathrm{GCD}_n$ must be exponentially equal to $q_n$, so that both sides of the inequality $\left|\frac{\prod_{i=1}^{n} b_i}{q_{n+1}q_n}\right| \leq \frac{1}{\left|\frac{q_n}{\mathrm{GCD}_n}\right|^{1+\delta}}$ will decrease sub-exponentially. This requires a more delicate analysis that is outside the scope of this paper. In conclusion, we conjecture the condition

$$\mathrm{GCD}_n \doteq q_n \doteq B^{\frac{n}{2}} \cdot n!^{\frac{d_b}{2}}$$

to be necessary but not sufficient for yielding a nontrivial $\delta$ in this case.

Case 2: $\frac{d_b}{d_a} = 2$, $q_n \doteq \alpha^n \cdot n!^{d_a}$

Note that

$$\left|\frac{\prod_{i=1}^{n} b_i}{q_{n+1}q_n}\right| \doteq \left|\frac{B^n \cdot n!^{d_b}}{(\alpha^n \cdot n!^{d_a})^2}\right| = \left|\frac{B}{\alpha^2}\right|^n,$$

and since $|B| < \alpha^2$, the convergence rate is exponential. For the expression for $\delta$, we have

$$\frac{\ln|q_{n+1}| - \ln|\prod_{i=1}^{n} b_i| + \ln \mathrm{GCD}_n}{\ln|q_n| - \ln \mathrm{GCD}_n} =$$

$$\frac{(n+1)\ln|\alpha| + \boldsymbol{d_a(n+1)\ln(n+1)} + \ln\left|E_{n+1}^q\right| - n\ln B - \boldsymbol{2d_a n\ln n} - \ln|E_n^b| + \boldsymbol{\ln \mathrm{GCD}_n}}{n\ln|\alpha| + \boldsymbol{d_a n\ln n} + \ln\left|E_n^q\right| - \boldsymbol{\ln \mathrm{GCD}_n}}.$$

Dominant terms that determine the behavior of the expression are highlighted in bold. Here comes the crucial part of the proof for **Theorem 1**: To provide a nontrivial $\delta$, the limsup of this expression must not be $-1$. Thus, $\mathrm{GCD}_n$ must contain a super-exponential factor at the size $n!^{d_a}$ – i.e., **factorial reduction**, proving **Theorem 1**. Note that a bigger super-exponential factor is not possible since $p_n \doteq q_n \doteq \alpha^n \cdot n!^{d_a}$.

Now, to prove **Theorem 2**, inserting $\mathrm{GCD}_n \doteq \lambda^n \cdot n!^{d_a}$ into the last expression yields

$$\frac{\boldsymbol{(n+1)\ln|\alpha|} + d_a(n+1)\ln(n+1) + \ln\left|E_{n+1}^q\right| - \boldsymbol{n\ln|B|} - 2d_a n\ln n - \ln|E_n^b| + \boldsymbol{n\ln\lambda} + d_a n\ln n + \ln|E_n^{\mathrm{GCD}}|}{\boldsymbol{n\ln|\alpha|} + d_a n\ln n + \ln\left|E_n^q\right| - \boldsymbol{n\ln\lambda} - d_a n\ln n - \ln|E_n^{\mathrm{GCD}}|}$$

$$\xrightarrow[n\to\infty]{} \frac{\ln|\alpha| - \ln|B| + \ln\lambda}{\ln|\alpha| - \ln\lambda}.$$

The terms containing multiples of $n\,ln\,n$ canceled out and the ones determining the behavior of this expression are highlighted in bold. This matches the expression for $\delta$ in **Theorem 2** and proves it.

<u>Case 3</u>: $\frac{d_b}{d_a} < 2$, $A^n \cdot n!^{d_a}$

Note that

$$\left|\frac{\prod_{i=1}^n b_i}{q_{n+1}q_n}\right| \doteq \left|\frac{B^n \cdot n!^{d_b}}{(A^n \cdot n!^{d_a})^2}\right| = \left|\frac{B}{A^2}\right|^n \cdot n!^{d_b - 2d_a},$$

and since $d_b - 2d_a < 0$, the convergence rate is super-exponential. For the expression for $\delta$, even if $\mathrm{GCD}_n = 1$, we have

$$\frac{\ln|q_{n+1}| - \ln|\prod_{i=1}^n b_i| + \ln\mathrm{GCD}_n}{\ln|q_n| - \ln\mathrm{GCD}_n} =$$

$$\frac{(n+1)\ln|A| + \boldsymbol{d_a(n+1)\ln(n+1)} + \ln\left|E_n^q\right| - n\ln B - \boldsymbol{d_b n\ln n} - \ln|E_n^b| + 0}{\boldsymbol{d_a n\ln n} + n|A| + \ln\left|E_n^q\right| - 0}$$

$$\xrightarrow[n\to\infty]{} \frac{d_a - d_b}{d_a},$$

which is positive, and proves irrationality, if and only if $d_a > d_b$.

For $\mathrm{GCD_n}$ to modify the limit above, it must be of factorial order. If so, and the factorial power is of size $n!^r$ (the exponential factors have no effect), then

$$\frac{\ln|q_{n+1}| - \ln|\prod_{i=1}^{n} b_i| + \ln \mathrm{GCD}_n}{\ln|q_n| - \ln \mathrm{GCD}_n} =$$

$$\frac{(n+1)\ln|A| + \boldsymbol{d_a(n+1)\ln(n+1)} + \ln\left|E_n^q\right| - n\ln B - \boldsymbol{d_b n \ln n} - \ln|E_n^b| + 0}{\boldsymbol{d_a n \ln n} + n|A| + \ln\left|E_n^q\right| - 0}$$

$$\xrightarrow[n\to\infty]{} \frac{d_a - d_b + r}{d_a - r},$$

which is better (larger) than before. In conclusion, it is not necessary or sufficient for $\mathrm{GCD}_n$ to be of factorial order to provide a nontrivial $\delta$ in this case.

## Appendix C: Inflation and deflation of continued fractions

In his paper, Apéry showed a linear recursion of depth 2 with rational function coefficients (ratio of two polynomial) and a related PCF is presented in [2]. The direct translation of the Apéry's recursion into a continued fraction has $a_n$ and $b_n$ as rational functions and not integer polynomials. However, they can be converted to a PCF form. To see the conversion, we multiply $a_n$ and $b_n$ by a non-zero sequence, thus converting them to integer polynomials without changing the limit or the number $p_n/q_n$ at each $n$. We call this process "inflation". This process is also needed when some rational coefficient is used in **Conjecture 1.4**.

Conversely, any PCF that has been multiplied by a non-zero sequence can be simplified by removing that sequence. We call this process "deflation". Deflating makes the PCFs' $b_n$ and $a_n$ smaller (possibly of a lower degree), and most importantly, helps simplify $\mathrm{GCD}_n$s, despite not changing the induced $\delta$. This process can explain some powers of prime in Table 1.

**Property 1** *(Inflation and deflation of continued fractions). Let $c_n$ be a sequence of non-zero numbers. For every $n$*

$$a_0+\cfrac{b_1}{a_1+\cfrac{b_2}{a_2+\cfrac{b_3}{a_3\,...+\cfrac{b_n}{a_n}}}}=a_0+\cfrac{c_1b_1}{c_1a_1+\cfrac{c_1c_2b_2}{c_2a_2+\cfrac{c_2c_3b_3}{c_3a_3\,...+\cfrac{c_{n-1}c_nb_n}{c_na_n}}}}.$$

*I.e., if $p_n$ and $q_n$ are the numerator and denominator of the left hand side and $p_n'$ and $q_n'$ are those of the right hand side, then for every $n$,*

$$\frac{p_n}{q_n} = \frac{p'_n}{q'_n}.$$

Example: the relation between Apéry’s recursion and his PCF

Setting $a_n = \frac{34n^3+51n+27n+5}{(n+1)^3}$ and $b_n = -\frac{n^3}{(n+1)^3}$, Apéry used the recursion

$$u_{n+1} = a_{n+1}u_n + b_{n+1}u_{n-1}$$

with the initial conditions

$$p_{-1} = 1, \;\; p_0 = 5$$

$$q_{-1} = 0, \;\; q_0 = 1,$$

which generates the following rational continued fraction:

$$\frac{6}{\zeta(3)} = 5 - \cfrac{\frac{1}{8}}{\frac{117}{8} - \cfrac{\frac{8}{27}}{\frac{535}{27} - \cfrac{\frac{27}{64}}{\frac{1463}{64} \cdots - \cfrac{\frac{n^3}{(n+1)^3}}{\frac{34n^3 + 51n + 27n + 5}{(n+1)^3}}}}}.$$

Applying **Property 1**, we can inflate this continued fraction using the denominators of $a_n$ and $b_n$, i.e., the sequence $c_n = (n+1)^3$, and obtain the PCF

$$\frac{6}{\zeta(3)} = 5 - \cfrac{\frac{1}{8}\cdot 2^3}{\frac{117}{8}\cdot 2^3 - \cfrac{\frac{8}{27}\cdot 2^3\cdot 3^3}{\frac{535}{27}\cdot 3^3 - \cfrac{\frac{27}{64}\cdot 4^3\cdot 3^3}{\frac{1463}{64}\cdot 4^3 \ldots - \cfrac{\frac{n^3}{(n+1)^3}\cdot n^3(n+1)^3}{\frac{34n^3+51n+27n+5}{(n+1)^3}\cdot (n+1)^3}}}}$$

$$= 5 - \cfrac{1}{117 - \cfrac{64}{535 - \cfrac{279}{1463 \ldots - \cfrac{n^6}{34n^3+51n^2+27n+5}}}},$$

that is, [2]'s PCF, which is presented in our introduction.

To prove **Property 1**, and examine effect of this process on $\mathrm{GCD}_n$, we present the following lemma.

**Lemma 5.** *Consider the recursions*

$$u'_n = a'_n u'_{n-1} + b'_n u'_{n-2}$$

$$u_n = a_n u_{n-1} + b_n u_{n-2},$$

*where* $a'_n = c_n \cdot a_n$ *and* $b'_n = c_{n-1}c_n \cdot b_n$ *for some non-zero sequence* $c_n$. *If the recursion have the same initial values, we have*

$$u'_n = \left(\prod_{i=1}^{n} c_i\right)\cdot u_n.$$

*Proof.* Initializing the induction at $n = -1,0,$ is trivial since the product is empty.

To prove the induction at $n+1,$ write

$$u'_{n+1} = a'_{n+1}u'_n + b'_{n+1}u'_{n-1} = c_{n+1}a_{n+1}u'_n + c_n c_{n+1}b_{n+1}u'_{n-1} \underbrace{=}_{\text{assumption}}$$

$$c_{n+1}a_{n+1}\left(\prod_{i=1}^{n} c_i\right)u_n + c_n c_{n+1}b_{n+1}\left(\prod_{i=1}^{n-1} c_i\right)u_{n-1} =$$

$$\left(\prod_{i=1}^{n+1} c_i\right)(a_{n+1}u_n + b_{n+1}u_{n-1}) = \left(\prod_{i=1}^{n+1} c_i\right)u_{n+1}. \qquad \square$$

As a corollary, we can prove **Property 1** since

$$\frac{p'_n}{q'_n} = \frac{(\prod_{i=1}^{n} c_i)p_n}{(\prod_{i=1}^{n} c_i)q_n} = \frac{p_n}{q_n}.$$

Also, we state the following property, showing effect of this process on $\mathrm{GCD}_n$.

**Property 2.** *Using **Lemma 5** notation, for an inflated* PCF*, or any inflated rational* GCF,

$$\mathrm{GCD}[p'_n, q'_n] = \mathrm{GCD}\left[\left(\prod_{i=1}^{n} c_i\right)p_n, \left(\prod_{i=1}^{n} c_i\right)q_n\right] = \left(\prod_{i=1}^{n} c_i\right)\mathrm{GCD}[p_n, q_n].$$

<u>Example</u>:

Considering the following PCF from Table 1:

$$\mathrm{PCF}[a'_n, b'_n] = 1 + \cfrac{3}{3 + \cfrac{15}{5 + \cfrac{35}{7 \ldots + \cfrac{(3n+1)(3n-2)}{3n+1}}}},$$

and observe that it is inflated by the sequence $c_n = 3n + 1$. By deflating it, we obtain the *regular* continued fraction of the golden ratio $\varphi$:

$$1 + \cfrac{1}{1 + \cfrac{1}{1 + \cfrac{1}{1 \ldots + \frac{1}{1}}}},$$

which upholds $p_n = F_{n+2}, q_n = F_{n+1}$ , where $F_n$ is the $n$th Fibonacci number. Therefore, for the original PCF we obtain

$$\mathrm{GCD}[p'_n, q'_n] = \left(\prod_{i=1}^{n} c_i\right) \mathrm{GCD}[p_n, q_n] = \prod_{i=1}^{n} (3i + 1) = (3n + 1)!!!$$

since consecutive Fibonacci numbers are coprime.

Inflation by "$\sqrt{p}$":

An interesting special case of deflation occurs when for some integer $r$, we have $r \mid a_n$, $r \mid b_n$, but $r^2 \nmid b_n$ (so we do not have a trivial inflation by $c_n = r$). In this case, $\mathrm{GCD}_n$ is divided by powers of $r$ of the forms $r^{\lfloor n/2 \rfloor}$. For instance, for PCF$[3n + 6, 3n^2 + 9n]$, we found (numerically)

$$\mathrm{GCD}_n \doteq n! \cdot \sqrt{3}^n$$

We call this phenomenon inflation with $\sqrt{p}$, and it can be explained by considering inflation with the following sequence

$$c_n = \begin{cases} r & if \quad n \; is \; even \\ 1 & if \quad n \; is \; odd \end{cases}.$$

## Appendix D: Analysis and proof of $\text{GCD}_n$ formula for some PCFs

This section concentrate on various PCFs for which we identified and proven an explicit part of $\text{GCD}_n$. As example 2 in section 2.2 shows, for certain families of PCFs this part prove the existence of FR with sufficiently large $\lambda$ for proving irrationality. This proof is presented here hoping it will promote future research of additional proofs.

**Theorem 5.** *If* $\text{PCF}[a_n, b_n]$ *satisfies for all* $n$

$$a_n = -a_{-1-n}, \tag{Eq. D1}$$

$$b_n = b_{-n}, \tag{Eq. D2}$$

$$b_0 = 0, \tag{Eq. D3}$$

*then for all* $n$

$$\frac{n!}{\text{LCM}[n]} \Bigg| \text{GCD}_n \cdot 2^n.$$

Note that conditions in $(\text{Eq. D1,2})$ involves negative-indexed coefficients that are not used (or defined) by the PCF. Nevertheless, since $a_n$ and $b_n$ are given by polynomials, we can extend them to all values of $n$.

*Proof*. Recall that both $p_n$ and $q_n$ satisfy the recursion

$$u_n = a_n u_{n-1} + b_n u_{n-2}. \tag{Eq. D4}$$

We first prove that for any odd prime $p$, $u_n$ is divisible by $p$ for all $n \geq 2p$. We do so by analyzing the sequence $u_n$ modulo $p$.

Let $p = 2h + 1$, we prove by induction that for all $-1 \leq m < h$

$$u_{h+m} = u_{h-m-2} \cdot b_h \cdot b_{h-1} \cdot \ldots \cdot b_{h-m} \pmod{p}. \qquad \text{(Eq. D5)}$$

Initializing the induction at $m = -1$ is trivial since $u_{h-1} = u_{h-1}$.

Initializing the induction at $m = 0$ requires $u_h = u_{h-2} b_h \pmod{p}$. To prove that we start by using (Eq. D1) at $n = h$, i.e.,

$$a_h = -a_{-1-h}.$$

Note that for modulo $r$, every polynomial has a period of $p = 2h + 1$, namely

$$-a_{-1-h} = -a_h = 0 \pmod{p}.$$

Substituting this into the recursion (Eq. D4) at $n = h$, we obtain as required

$$u_h = a_h u_{h-1} + b_h u_{h-2} = b_h u_{h-2} \pmod{p}.$$

To prove the induction at $m + 1$, write (Eq. D1,2) for $n = h + m + 1$

$$a_{h+m+1} + a_{-h-m-2} = 0,$$

$$b_{h+m+1} - b_{-h-m-1} = 0.$$

Using periodicity modulo $p = 2h + 1$,

$$a_{-h-m-2} = a_{h-m-1} \pmod{p},$$

$$b_{-h-m-1} = b_{h-m} \pmod{p},$$

and combining this with the symmetries, we obtain

$$a_{h+m+1} = -a_{h-m-1} \pmod{p},$$

$$b_{h+m+1} = b_{h-m} \pmod{p}.$$

Now, substitute these into the recursion (Eq. D4) at $n = h + m + 1$ gets

$$u_{h+m+1} = a_{h+m+1} u_{h+m} + b_{h+m+1} u_{h+m-1} = -a_{h-m-1} u_{h+m} + b_{h-m} u_{h+m-1} \pmod p.$$

Using the induction assumption at $m$ and $m - 1$, we have

$$u_{h+m+1} = -a_{h-m-1} u_{h-m-2} b_h b_{h-1} \dots b_{h-m} + b_{h-m} u_{h-m-1} b_h b_{h-1} \dots b_{h-m+1} \pmod p,$$

and by rearranging, we get

$$u_{h+m+1} = (-a_{h-m-1} u_{h-m-2} + u_{h-m-1}) b_h b_{h-1} \dots b_{h-m} \pmod p.$$

Substituting the recursion (Eq. D4) at $n = h - m - 1$ achieve the induction at $m + 1$

$$u_{h+m+1} = u_{h-m-3} b_h b_{h-1} \dots b_{h-m-1} \pmod p,$$

completing the proof of the induction.

Using the periodicity of $a_n$ and $b_n$ modulo $p$, the relation

$$u_{h+m} = u_{h-m-2} \cdot b_h \cdot b_{h-1} \cdot \dots \cdot b_{h-m} \pmod r$$

can be proven also after shifting by $p$, i.e.,

$$u_{p+h+m} = u_{p+h-m-2} \cdot b_{p+h} \cdot b_{p+h-1} \cdot \dots \cdot b_{p+h-m} \pmod p.$$

Now, set $h \leq m < p + h$. Then $b_p$ exists among the multipliers of the right hand side. Combining periodicity and (Eq. D3), implying $b_p = b_0 = 0 \pmod p$, we finally have $u_{p+h+m} = 0 \pmod p$. So in other words, $p$ divides $u_n$ starting at $n = p + 2h = 2r - 1$.

We now explain why this result suffices to prove the theorem. We proved that an odd prime $p$ divides both $p_n$ and $q_n$ for $n \geq 2p - 1$, and therefore, it divides $\mathrm{GCD}_n$, in line with the fact that $p$ divides $\frac{n!}{\mathrm{LCM}[n]}$ for $n \geq 2p$. Furthermore, since $u_n$ is divisible by $p$ for $n \geq 2p-1$, we can set a new sequence

$$u'_n = \frac{u_{n+2p}}{p},$$

which is well defined for $n \geq -1$. Since the sequence is obtained by shifting (by a multiple of $p$) and scaling $u_n$, it satisfies the same recursion $(\mathrm{D}.4)$ at $n \geq 1$. However, unlike the original sequence, $u'_n$ in fact also satisfies the recursion also at $n = 0$. This follows from the fact that $b_0 = 0$, and therefore the recursion at $n = 0$ does not involve $u'_{-2}$. Note that the original sequence $u_n$ has an arbitrary initial condition at $u_0$ and $u_{-1}$ that may not satisfy the recursion. This is not the case for $u'_n$ since it was generated by the recursion even at $n = 0$.

Following that, we apply the above induction result $(\mathrm{Eq.\,D5})$ and prove that $p$ divides $u'_n$ starting at $n = p - 1$ (instead of $n = 2p - 1$ as in the original $u_n$). In other words, $p^2$ divides the original $u_n$ starting at $n = 3p - 1$. In general, $p^k$ divides $u_n$ starting at $n = (k+1)p - 1$. Note that for $p > \sqrt{n}$, $p^k$ divides $\frac{n!}{2^l \mathrm{LCM}[n]}$ starting at $n = (k+1)p$, and therefore, the requirement is met. Moreover, $n!$ obtains an additional $p$ factor at $n = k \cdot p^2$. For $k = 1$, this factor is canceled by the denominator's $\mathrm{LCM}[n]$, but for $k > 1$ we must prove this additional factor of $\mathrm{GCD}_n$.

To do that, observe that the original proof that $p$ divides $u_n$ for $n \geq 2p - 1$ was valid also when substituting $p$ with any power of prime. Thus, at $n = 2p^2 - 1$, the sequence $u_n$ in fact obtain a factor of $p^2$ instead of the expected single factor of $p$. Appling this for any prime power fulfills the requirements of the theorem.

## Appendix E: Additional examples of PCFs with factorial reduction

One may wonder whether the conjectures discovered in this study are indeed mathematical truth or merely mathematical coincidences that break down at higher degrees or larger coefficients. However, the method employed in this study makes it fairly unlikely that the conjectures will break down. Nevertheless, such an assumption does not replace the need for formal proof. We believe that many (if not all) of the new conjectures are indeed truths awaiting rigorous proof, relying on vast search spaces examined in this work. To strengthen our conjectures, we give additional examples abundantly. Moreover, this appendix addresses three further causes:

1. Visualizing the results.
2. Present a piece of the yet unknown general structure toward revealing it whole.
3. Lay the groundwork and provide more data for proofs or additional conjectures.

<u>FR and rational roots</u>

**Conjecture 1.1** states that, for a given $b_n$, there exists an $a_n$ such that $\mathrm{PCF}[a_n, b_n]$ has FR if and only if $b_n$ has only rational roots. **Table 3** shows a classification, by numerical tests, of *all* the $b_n$ polynomials of degree 2 and with integer coefficients between 1 and 4. For each such $b_n$, we search for $a_n$ polynomials in the integer coefficient range 1 to 5.

| Have FR | | | Do not have FR | | |
|---|---|---|---|---|---|
| $b_n$ | $b_n$'s roots | $a_n$ example | $b_n$ | $b_n$'s roots | $a_n$ example |
| $n^2 + 2n + 1$ | $-1$<br>$-1$ | $3 + 2n$ | $n^2 + 4n + 2$ | $-2 - \sqrt{2}$<br>$-2 + \sqrt{2}$ | none<br>found |
| $2n^2 + 3n + 1$ | $-1$<br>$-1/2$ | $n + 1$ | $n^2 + 4n + 1$ | $-2 - \sqrt{3}$<br>$-2 + \sqrt{3}$ | none<br>found |
| $3n^2 + 4n + 1$ | $-1$<br>$-1/3$ | $1 + 2$ | $n^2 + 1$ | $i$<br>$-i$ | none<br>found |
| $4n^2 + 4n + 1$ | $-1/2$<br>$-1/2$ | $n + 1$ | $n^2 + 2n + 2$ | $-1 - i$<br>$-1 + i$ | none<br>found |
| $n^2 + 3n + 2$ | $-2$<br>$-1$ | $2n + 4$* | $4n^2 + 4n + 2$ | $-1/2 - i/2$<br>$-1/2 + i/2$ | none<br>found |
| $2n^2 + 4n + 2$ | $-1$<br>$-1$ | $n + 1$ | $n^2 + 2n + 3$ | $-1 - i\sqrt{2}$<br>$-1 + i\sqrt{2}$ | none<br>found |
| $n^2 + 4n + 3$ | $-3$<br>$-1$ | $5 + 2n$ | everything else | not rational | none<br>found |
| $n^2 + 4n + 4$ | $-2$<br>$-2$ | $2n + 5$ | … | … | … |

**Table 3: $b_n$ from degree 2 and coefficients 1 to 4, classified by the existence of FR.** We can see that, for degree 2, only the $b_n$ polynomials with FR have all-rational roots and vice versa. Further, note that all the $a_n$ examples belong to the conjectured complete structure for $\deg b_n = 2, \deg a_n = 1$ (**Conjecture 1.4**). *This PCF is an inflation of the *regular* continued fraction of $\sqrt{2} + 1$.

A piece of structure for all degrees

In this section we urge a generalizing of **Conjecture 1.4**, that deals only with $b_n$s of degree 2, by demonstrating families of balanced-degree PCFs with $b_n$s of degrees$>$ 2 and that have FR.

| $d_b = 2$ | $d_b = 4$ | $d_b = 4$ | $d_b = 4$ | $d_b = 6$ | $d_b = 6$ |
|---|---|---|---|---|---|
| $b_n = -n^2$ | $b_n = -n^4$ | $b_n = -4n^4 - 6n^3$ | $b_n = -n^2 \cdot (n+2)(2n-3)$ | $b_n = 4n^6 - 2n^5$ | $-n^6$ |
| $a_n$'s | $a_n$'s | $a_n$'s | $a_n$'s | $a_n$'s | $a_n$'s |
| $2n+1$ | $2n^2+2n+1$ | $4n^2+7n+2$ | $3n^2+n+3$ | $3n^3-5n^2-3n-1$ | $2n^3+3n^2+3n+1$ |
| $4n+2$ | $2n^2+2n+3$ | $4n^2+7n+3$ | $3n^2+3n-5$ | $3n^3+10n^2+8n+2$ | $2n^3+3n^2+11n+5$ |
| $6n+3$ | $2n^2+2n+7$ | $4n^2+7n+5$ | $3n^2+3n-1$ | | $6n^3+9n^2+5n+1$ |
| $8n+4$ | $2n^2+2n+13$ | $4n^2+7n+8$ | $3n^2+5n-3$ | | |
| $10n+5$ | | $4n^2+7n+12$ | $3n^2+7n+3$ | | |
| $12n+6$ | | $4n^2+7n+17$ | $3n^2+15n+13$ | | |
| $14n+7$ | | $5n^2+8n+1$ | $3n^2+15n+17$ | | |
| $16n+8$ | | $5n^2+14n+10$ | | | |

**Table 4: $b_n$ polynomials and (numerically) found $a_n$s for which PCF$[a_n, b_n]$ has FR.** For the first example, where $\deg(b_n) = d_b = 2$, we conjecture a structure in **Conjecture 1.4**. However, we still do not possess a solid conjecture for higher degrees.

More rational recurrence relations that yield integer sequences

As follows from Section 2.6, a conjectured $\mathrm{GCD}_n$ formula can be proven to be true if the following recurrence yields coprime integer sequences for any linearly-independent pairs of integer initial values:

$$u'_n = \frac{\mathrm{GCD}_{n-1}}{\mathrm{GCD}_n} a_n u'_{n-1} + \frac{\mathrm{GCD}_{n-2}}{\mathrm{GCD}_n} b_n u'_{n-2}.$$

We address the community to prove this property for simple $\mathrm{GCD}_n$s with the following examples. Furthermore, we request general conditions on $a_n$ and $b_n$ so that this property holds, either for the presented special cases or hopefully for other cases (of **Conjecture 2**).

| $\mathrm{GCD}_n \doteq n!$ <br> [1] $n(n-1)u'_n = (n-1)a_n u'_{n-1} + b_n u'_{n-2}$ | | | | $\mathrm{GCD}_n \doteq n!/2^n$ <br> [1] $n(n-1)u'_n = 2(n-1)a_n u'_{n-1} + 2b_n u'_{n-2}$ | | | |
|---|---|---|---|---|---|---|---|
| $a_n$ | $b_n$ | $a_n$ | $b_n$ | $a_n$ | $b_n$ | $a_n$ | $b_n$ |
| $3n+1$ | $4n^2+2n$ | $n+3$ | $2n^2+4n$ | $n-2$ | $2n^2-n$ | $3n+1$ | $-2n^2-1$ |
| $-n+1$ | $6n^2+12n$ | $2n+4$ | $3n^2+9n$ | $3n-2$ | $4n^2-n$ | $-n+2$ | $2n^2-n$ |
| $n-1$ | $6n^2+12n$ | $7n+3$ | $8n^2+4n$ | $-n-2$ | $6n^2+3n$ | $5n+1$ | $6n^2+9n$ |
| $7n+1$ | $8n^2+2n$ | $8n+2$ | $9n^2+3n$ | $7n-2$ | $8n^2-n$ | $-2n+1$ | $8n^2+2n$ |
| $3n+1$ | $10n^2+20n$ | $15n+1$ | $16n^2+2n$ | $2n-1$ | $8n^2+2n$ | $4n+1$ | $12n^2+6n$ |
| $-n+1$ | $12n^2+6n$ | $15n+3$ | $16n^2+4n$ | $n-2$ | $12n^2+3n$ | $-n+2$ | $12n^2+3n$ |
| $n-1$ | $12n^2+6n$ | $17n+5$ | $18n^2+6n$ | $15n-2$ | $16n^2-n$ | $6n-1$ | $16n^2+2$ |

**Table 5: $a_n$ and $b_n$ examples for special simplified recurrence relations yielding integer sequences.** A proof of this property also proves the formula for $\mathrm{GCD}_n$. Note that these are special cases as Section 2.6, and can be generalized to all $\mathrm{GCD}_n$ forms. [1]If the expressions used are only exponentially equal to $\mathrm{GCD}_n$, and not exactly equal, the simplified recursions formulas might generate rational sequences with sub-exponential denominators.